\documentclass[10pt] {amsart}

\usepackage{geometry}
\usepackage{amssymb}
\usepackage{amsmath}
\usepackage{amsthm}
\usepackage{amsfonts}
\usepackage{tikz}
\usetikzlibrary{matrix}
\usepackage[latin1]{inputenc}

\newcommand\restr[2]{{
  \left.\kern-\nulldelimiterspace 
  #1 
  \vphantom{\big|} 
  \right|_{#2} 
  }}
\newcommand{\suchthat}{\;\ifnum\currentgrouptype=16 \middle\fi|\;}

\newtheorem{theorem}{Theorem}[section]

\newtheorem{proposition}[theorem]{Proposition}
\newtheorem{corollary}[theorem]{Corollary}

\theoremstyle{definition}
\newtheorem{definition}[theorem]{Definition}

\begin{document}
\def\hii{3}
\title{Dirac Surfaces and Threefolds}

\author{Geoffrey Scott}
\email{gscott@math.utoronto.edu}

\begin{abstract}
We describe Dirac structures on surfaces and 3-manifolds. Every Dirac structure on a surface $M$ is described either by a regular 1-foliation or by a section of a circle bundle obtained as a fiberwise compactification of the line bundle $\wedge^2TM$. Every Dirac structure on a 3-manifold $M$ is either the union of a presymplectic manifold and a foliated Poisson manifold, or the union of a Poisson manifold and a foliated presymplectic manifold.
\end{abstract}

\maketitle

\section{Introduction}

A 2-form $\omega$ on a manifold $M$ induces a skew-symmetric map $TM \rightarrow T^*M$ whose graph is an $n$-dimensional subbundle of $\mathbb{T}M := TM \oplus T^*M$, isotropic with respect to the natural inner product on $\mathbb{T}M$. Likewise, a bivector $\Pi$ on $M$ induces a skew-symmetric map $T^* \rightarrow T$ whose graph is also an $n$-dimensional isotropic subbundle of $\mathbb{T}M$. The language of Dirac geometry allows us to simultaneously generalize presymplectic geometry and Poisson geometry by taking our main geometric structure to be an isotropic $n$-dimensional subspace of $\mathbb{T}M$, involutive with respect to the Courant bracket on $\mathbb{T}M$. The involutivity criterion generalizes the condition of a 2-form being closed, and of a bivector being Poisson. We review the language of Dirac geometry in Section \ref{sec:diracgeometry}. In Section \ref{sec:dim2}, we classify Dirac structures on surfaces. In Section \ref{sec:dim3}, we classify Dirac structures on 3-manifolds.

\section{Dirac Geometry}\label{sec:diracgeometry}
Proofs of the assertions made in this section can be found in \cite{C}, \cite{CW}, or any introduction to Dirac geometry.

\subsection{Linear Algebra}
For a real vector space $V$ of dimension $n$, $V \oplus V^*$ has a natural split-signature inner product $\langle v + \xi, w + \eta \rangle = \frac{1}{2} (\xi(w) + \eta(v))$ and natural projections $\rho: V \oplus V^* \rightarrow V$ and $\widehat{\rho}: V \oplus V^* \rightarrow V^*$. A subspace $L$ of $V \oplus V^*$ is \textbf{Lagrangian} if it is $n$-dimensional and isotropic with respect to this inner product, and $\textrm{Lag}(V \oplus V^*)$ is defined to be the space of all Lagrangians in $V \oplus V^*$. For $L \in \textrm{Lag}(V \oplus V^*)$, let $\Delta := \rho(L)$ and $\widehat{\Delta} := \widehat{\rho}(L)$. The flags of subspaces in $V$ and $V^*$
\[
0 \subseteq L \cap V \subseteq \Delta \subseteq V \hspace{1cm} \textrm{and} \hspace{1cm} 0 \subseteq L \cap V^* \subseteq \widehat{\Delta} \subseteq V^*
\]
are related by the fact that $\textrm{Ann}(L \cap V) = \widehat{\Delta}$ and $\textrm{Ann}(\Delta) = L \cap V^*$. The \textbf{type} of $L \in \textrm{Lag}(V \oplus V^*)$ is the pair of integers $(\textrm{dim}(L \cap V), \textrm{dim}(L \cap V^*))$. For example, the type of $V$ is $(n, 0)$, and the type of a Lagrangian defined as the graph of a 2-form $\omega \in \wedge^2V^*$ is $(\textrm{dim}(\textrm{ker} \ \omega ), 0)$. In most references, the \emph{type} of $L$ is instead defined as the codimension of $\Delta$ in $V$, which equals the second coordinate of our definition of type. Our definition distinguishes between $L = V \subseteq V \oplus V^*$ and the graph of a symplectic structure, but our definition is not well-defined on general Courant algebroids that do not contain $TM$ as a subbundle. We call a Lagrangian of type $(a, b)$ \textbf{even} if $b$ is even, and \textbf{odd} otherwise. Topologically, $\textrm{Lag}(V \oplus V^*) \cong O(n)$ (\cite{C}, Section 1.3); its two components are the even Lagrangians $\textrm{Lag}_e(V \oplus V^*) \subseteq \textrm{Lag}(V \oplus V^*)$ and the odd Lagrangians $\textrm{Lag}_o(V \oplus V^*) \subseteq \textrm{Lag}(V \oplus V^*)$. A Lagrangian $L$ defines a skew form $\epsilon \in \wedge^2 \Delta^*$ and a bivector $\Pi \in \wedge^2 \widehat{\Delta}^* \cong \wedge^2 (V / (L \cap V))$ by the formulas
\[
\epsilon(\rho(x)) = \restr{\widehat{\rho}(x)}{\Delta}
\hspace{0.75cm} \textrm{and} \hspace{0.75cm}
\Pi(\widehat{\rho}(x)) = \restr{\rho(x)}{\widehat{\Delta}} \hspace{1cm} \textrm{for all $x \in V \oplus V^*$}
\]
Conversely, any Lagrangian is uniquely specified by the pair $(\Delta, \epsilon)$, and also by the pair $(\widehat{\Delta}, \Pi)$.
\begin{proposition}\label{prop:cour} (\cite{C}, Proposition 1.1.5) The maps $L \mapsto (\Delta, \epsilon)$ and $L \mapsto (\Delta^*, \Pi)$ define bijections
\[
\left\{ \begin{array}{c} \textrm{Pairs} \ (\Delta, \epsilon)\\ \Delta \ \textrm{a subspace of} \ V \\ \epsilon \in \wedge^2\Delta^* \end{array} \right\} \leftrightarrow \textrm{Lag}(  V \oplus V^*) \leftrightarrow \left\{ \begin{array}{c} \textrm{Pairs} \ (\widehat{\Delta}, \Pi)\\ \widehat{\Delta} \ \textrm{a subspace of} \ V^* \\ \Pi \in \wedge^2\widehat{\Delta}^* \end{array} \right\}
\]
\end{proposition}

\subsection{Differential Geometry}
The \textbf{generalized tangent bundle} of a manifold $M$ is the bundle $\mathbb{T}M := TM \oplus T^*M$ endowed with the split-signature inner product $\langle X + \xi, Y + \eta \rangle = \frac{1}{2}(\eta(X) + \xi(Y))$ and the \textbf{Courant bracket}
\[
[X + \xi, Y + \eta] = [X, Y] + \mathcal{L}_X\eta - \mathcal{L}_Y\xi - \frac{1}{2}d(\iota_X\eta - \iota_Y\xi).
\]
on its space of sections. We denote by $\textrm{Lag}(\mathbb{T}M)$ the $O(n)$ bundle over $M$ whose fiber over $x \in M$ is $\textrm{Lag}(T_xM \oplus T_x^*M)$. An \textbf{almost Dirac structure} is a Lagrangian subbundle $E$ of $\mathbb{T}M$ -- equivalently, a section of $\textrm{Lag}(\mathbb{T}M)$. Given an almost Dirac structure $E$ on $M$, we denote
\[
M_{(a, b)} := \{x \in M \mid E_x \ \textrm{has type} \ (a, b)\}.
\]
An almost Dirac struture $E$ is a \textbf{Dirac structure} if it is involutive -- that is, if the Courant bracket of two sections of $E$ is again a section of $E$. This involutivity criterion can be expressed using the language of foliated forms (or foliated poisson structures) on regions of $M$ where the dimension of $\Delta$ (or $\widehat{\Delta}$) is locally constant. 
\begin{definition}
Let $D$ be a distribution on a manifold $M$ given by a regular foliation. The \textbf{foliated exterior derivative} 
\[
d_D: \Gamma(\wedge^k D^*) \rightarrow \Gamma(\wedge^{k+1} D^*)
\]
is given by the usual Cartan formula
\begin{align*}
d_D\omega(V_0, \dots, V_k) &:= \sum_{i}(-1)^ia(V_i)\left( \omega(V_0, \dots, \hat{V}_i, \dots, V_k)\right)\\
& \ \ \ \  + \sum_{i < j} (-1)^{i+j} \omega([V_i, V_j], V_0, \dots, \hat{V}_i, \dots, \hat{V}_j, \dots, V_k).
\end{align*}
where $\{V_0, \dots, V_k\}$ are sections of $D$. A \textbf{foliated presymplectic} form is a $d_D$-closed $\omega \in \Gamma(\wedge^2 D^*)$. If $D = TM$, we recover the usual definition of a presymplectic form. 
\end{definition}
\begin{definition} Let $D$ be a distribution on a manifold $M$ given by a regular foliation. The sheaf of \textbf{admissible} functions on $M$ is the sheaf $C^{\infty}_D$ of functions which are constant on the leaves of $D$. A \textbf{foliated Poisson} structure is a Poisson bracket on $C^{\infty}_D$.
\end{definition}

\begin{proposition}\label{prop:integrability}(\cite{G} Prop. 2.7, \cite{C} Prop 2.5.3 and  Cor. 2.6.3) Let $E$ be an almost Dirac structure
\begin{enumerate}
\item If $\Delta$ is a subbundle of $TM$ (so $E$ can be described as the graph of $\epsilon \in \wedge^2(\Delta)^*$), then $E$ is a Dirac structure if and only if $\Delta$ integrates to a foliation and $d_{\Delta} \epsilon = 0$.
\item If $E \cap TM$ is a subbundle of $TM$ (so $E$ can be described as the graph of $\Pi \in \wedge^2(\widehat{\Delta})^*$), then $E$ is a Dirac structure if and only if $E \cap TM$ integrates to a foliation and $\Pi$ defines a foliated Poisson structure with on $M$ with respect to this foliation.
\end{enumerate}
\end{proposition}

For example, wherever an almost Dirac structure $E$ is transverse to $T^*M \subseteq \mathbb{T}M$, $E$ is the graph of a 2-form $\omega$ and Propostion \ref{prop:integrability} states that $E$ is Dirac structure precisely if $\omega$ is closed. Wherever $E$ is transverse to $TM$, $E$ is the graph of a bivector $\Pi$, and $E$ is Dirac precisely if $\Pi$ is Poisson.

\section{Dirac Structures on Surfaces}\label{sec:dim2}

\subsection{Linear Algebra} For a 2-dimensional vector space $V$, $\textrm{Lag}(V \oplus V^*) \cong O(2)$ consists of two circles, $\textrm{Lag}_e(V \oplus V^*)$ and $\textrm{Lag}_o(V \oplus V^*)$. The circle of even Lagrangians is covered by the maps
\begin{align}
\wedge^2V^* &\rightarrow \textrm{Lag}_e(V \oplus V^*)  \hspace{2cm} \textrm{and} &\wedge^2V &\rightarrow \textrm{Lag}_e(V \oplus V^*) \label{eqn:coverings}\\
\omega &\mapsto \textrm{Graph}(\omega)           &\Pi &\mapsto \textrm{Graph}(\Pi) \nonumber
\end{align}
Each even Lagrangian has type (0,0) except for $V$ and $V^*$, which are the images of $0 \in \wedge^2V^*$ and $0 \in \wedge^2V$ in the maps above, and have type $(2, 0)$ and $(0, 2)$, respectively. The circle of odd Lagrangians is isomorphic to $P(V)$, the projective space of lines in $V$, by
\begin{align*}
P(V) &\rightarrow \textrm{Lag}_o(V \oplus V^*)\\
L &\mapsto L + \textrm{Ann}(L)
\end{align*}

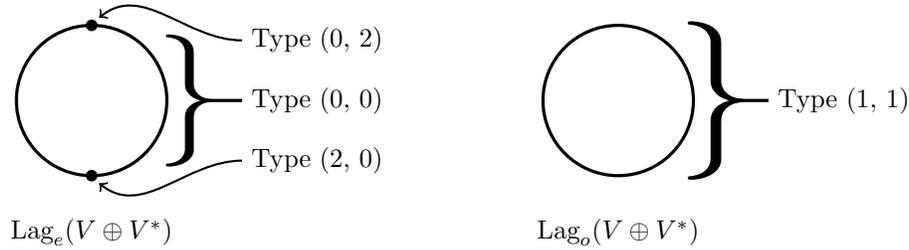
\begin{figure}[!ht]
\centering
\begin{tikzpicture}[scale = 1]
\begin{scope}[shift={(-3.5,0)}]
	\draw[very thick] (0, 0) circle(1cm);
	\draw[fill = black] (0, 1) circle(0.07cm);
	\draw[fill = black] (0, -1) circle(0.07cm);

  \node[anchor=west] at (2,0.8) (topdescription) {Type (0, 2)};
	\node[anchor=west] at (2,-0.8) (bottomdescription) {Type (2, 0)};
	\node[anchor=west] at (2, 0) (middledescription) {Type (0, 0)};
		
  \draw (topdescription) edge[thick, out=180,in=40,->] (0.1, 1.1);
	\draw (bottomdescription) edge[thick, out=180,in=-40,->] (0.1, -1.1);
	
	\node[scale = 5] at (1.3, 0) {$\}$};
	\draw (middledescription) edge[very thick] (1.5, 0);
	
	\node at (0, -1.75) {$\textrm{Lag}_e(V \oplus V^*)$};
		
\end{scope}

\begin{scope}[shift={(3.5,0)}]

	\draw[very thick] (0, 0) circle(1cm);

	\node[scale = 6] at (1.3, 0) {$\}$};
	\node[anchor=west] at (2, 0) (middledescription) {Type (1, 1)};
	\draw (middledescription) edge[very thick] (1.5, 0);
	
	\node at (0, -1.75) {$\textrm{Lag}_o(V \oplus V^*)$};
	
\end{scope}

\end{tikzpicture}
\caption{The topology of $\textrm{Lag}(V \oplus V^*)$ for $\textrm{dim}(V) = 2$}
\label{fig_disconnect}
\end{figure}

\subsection{Differential Geometry} Let $M$ be a surface. Then both $\textrm{Lag}_e(\mathbb{T}M)$ and $\textrm{Lag}_o(\mathbb{T}M)$ are circle bundles over $M$. The Maps (\ref{eqn:coverings}) globalize to inclusions
\[
\wedge^2TM \rightarrow \textrm{Lag}_e(\mathbb{T}M) \hspace{2cm} \textrm{and} \hspace{2cm} \wedge^2 T^*M \rightarrow \textrm{Lag}_e(\mathbb{T}M)
\]
so the circle bundle $\textrm{Lag}_e(\mathbb{T}M)$ may be viewed as the fiberwise compactification of the canonical (or anticanonical) line bundle.

Let $E$ be an even almost Dirac structure. On $M_{(0, 0)} \cup M_{(2, 0)}$, $E$ is the graph of a 2-form and Proposition \ref{prop:integrability} states that $E$ is Dirac if and only if this 2-form is closed. Every 2-form on a surface is closed, so $E$ is Dirac on $M_{(0, 0)} \cup M_{(2, 0)}$. 
Similarly, $E$ is Dirac on $M_{(0, 0)} \cup M_{(0, 2)}$ because every bivector on a surface is Poisson. Therefore, every even almost Dirac structure on a surface is Dirac. The data of an odd Dirac structure is equivalent to the data of a regular 1-dimensonial foliation on $M$. This proves
\begin{theorem}
Let $M$ be a surface. Even Dirac structures on $M$ are sections of the circle bundle $\textrm{Lag}_e(\mathbb{T}M)$. Odd Dirac structures on $M$ correspond to to regular $1$-foliations on $M$.
\end{theorem}
If $M$ is orientable, $\textrm{Lag}_e(\mathbb{T}M)$ is trivial, so $\Gamma(M, \textrm{Lag}_e(\mathbb{T}M)) \cong \textrm{Map}(M, S^1)$.
\begin{corollary}
Let $M$ be an orientable surface. The path components of the space of even Dirac structures are classified by $H^1(M; \mathbb{Z})$.
\end{corollary}
A necessary and sufficient condition for the existence of a regular 1-foliation on a surface is that the Euler characteristic of the surface vanishes.
\begin{corollary}
The only closed surfaces that admit an odd Dirac structure are the torus and the klein bottle.
\end{corollary}

\section{Dirac Structures on 3-manifolds}\label{sec:dim3}

\subsection{Linear Algebra} Let $V$ be a 3-dimensional real vector space, and let $Gr(k, V)$ denote the Grassmannian of $k$-dimensional planes in $V$. Then $\textrm{Lag}(V \oplus V^*) \cong O(3)$ is diffeomorphic to two copies of $\mathbb{R}P^3$. $\textrm{Lag}_e(V \oplus V^*)$ consists of Lagrangians of type $(1, 0), (3, 0)$, and $(1, 2)$. The Lagrangians of type $(1, 0)$ and $(1, 2)$ are precisely the ones for which $\textrm{dim}(\widehat{\Delta}) = 2$. By Proposition \ref{prop:cour}, these Lagrangians are classified by pairs $(\widehat{\Delta}, \Pi)$, where $\widehat{\Delta} \in \textrm{Gr}(2, V^*)$ and $\Pi \in \wedge^2\widehat{\Delta}^*$. The space of all such pairs is the total space of a real line bundle over $Gr(2, V^*) \cong \mathbb{R}P^2$ with projection $(\widehat{\Delta}, \Pi) \mapsto \widehat{\Delta}$. The zero section of this bundle corresponds to Lagrangians of type $(1, 2)$. Similarly, the Lagrangians of type $(1, 0)$ and $(3, 0)$ are precisely the Lagrangians for which $\Delta = V$, and are classified by elements $\epsilon \in \wedge^2V$. This set is a 3-ball whose zero corresponds to the Lagrangian $V$.

The space $\textrm{Lag}_o(V \oplus V^*)$ consists of Lagrangians of type $(0, 1), (0, 3)$, and $(2, 1)$. The Lagrangians of type $(0, 1)$ and $(2, 1)$ are classified by pairs $(\Delta, \epsilon)$, where $\Delta \in \textrm{Gr}(2, V)$ and $\epsilon \in \wedge^2\Delta^*$. This set is the total space of a real line bundle over $Gr(2, V) \cong \mathbb{R}P^2$ whose zero section corresponds to Lagrangians of type $(2, 1)$. The Lagrangians of type $(0, 1)$ and $(0, 3)$ are classified by elements $\epsilon \in \wedge^2V^*$, a 3-ball whose zero corresponds to $V^*$.

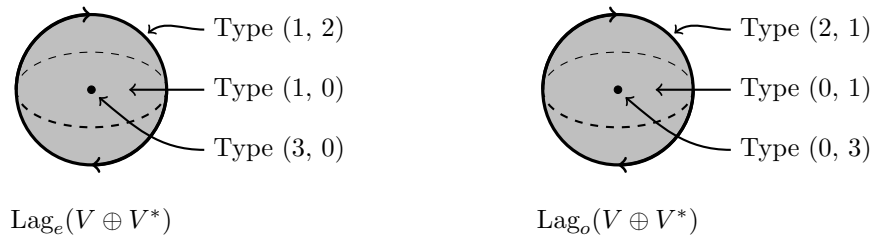
\begin{figure}[!ht]
\centering
\begin{tikzpicture}[scale = 1]

\begin{scope}[shift={(-3.5,0)}]
	\draw[very thick, fill = lightgray] (0, 0) circle(1cm);
	
	\draw[very thick, ->] (-1, 0) arc(180:90:1cm);
	\draw[very thick, ->] (1, 0) arc(0:-90:1cm);
	
	\draw[thick, dashed] (1, 0) arc (0:-180:1cm and 0.5cm);
	\draw[dashed] (1, 0) arc (0:180:1cm and 0.5cm);
	\draw[fill = black] (0, 0) circle(0.05cm);

  \node[anchor=west] at (1.5,0.8) (topdescription) {Type (1, 2)};
	\node[anchor=west] at (1.5,-0.8) (bottomdescription) {Type (3, 0)};
	\node[anchor=west] at (1.5, 0) (middledescription) {Type (1, 0)};
  \draw (topdescription) edge[thick, out=180,in=40,->] (0.75, 0.75);
	\draw (bottomdescription) edge[thick, out=180,in=-40,->] (0.1, -0.1);
	\draw (middledescription) edge[thick, ->] (0.5, 0);
	
	\node at (0, -1.75) {$\textrm{Lag}_e(V \oplus V^*)$};
\end{scope}

\begin{scope}[shift={(3.5,0)}]
	\draw[very thick, fill = lightgray] (0, 0) circle(1cm);
	
	\draw[very thick, ->] (-1, 0) arc(180:90:1cm);
	\draw[very thick, ->] (1, 0) arc(0:-90:1cm);
	
	\draw[thick, dashed] (1, 0) arc (0:-180:1cm and 0.5cm);
	\draw[dashed] (1, 0) arc (0:180:1cm and 0.5cm);
	\draw[fill = black] (0, 0) circle(0.05cm);

  \node[anchor=west] at (1.5,0.8) (topdescription) {Type (2, 1)};
	\node[anchor=west] at (1.5,-0.8) (bottomdescription) {Type (0, 3)};
	\node[anchor=west] at (1.5, 0) (middledescription) {Type (0, 1)};
  \draw (topdescription) edge[thick, out=180,in=40,->] (0.75, 0.75);
	\draw (bottomdescription) edge[thick, out=180,in=-40,->] (0.1, -0.1);
	\draw (middledescription) edge[thick, ->] (0.5, 0);
	
	\node at (0, -1.75) {$\textrm{Lag}_o(V \oplus V^*)$};
\end{scope}

\end{tikzpicture}
\caption{The topology of $\textrm{Lag}(V \oplus V^*)$ for $\textrm{dim}(V) = 3$.}
\label{fig:dim3_lag}
\end{figure}
\subsection{Differential Geometry} Let $E$ be a Dirac structure on a 3-manifold. We can use Proposition \ref{prop:integrability} and the discussion above to describe the Dirac geometry of different regions of $M$ using the language of foliated presymplectic and foliated Poisson structures
\begin{align*}
M_{(1, 0)} \cup M_{(1, 2)}&: \ \textrm{Foliated Poisson} & & 
M_{(1, 0)} \cup M_{(3, 0)}: \ \textrm{Presymplectic}\\
M_{(0, 1)} \cup M_{(2, 1)}&: \ \textrm{Foliated Presymplectic} & & 
M_{(0, 1)} \cup M_{(0, 3)}: \ \textrm{Poisson}
\end{align*}
This is summarized in the following theorem.
\begin{theorem} Let $M$ be a three dimensional manifold.
\begin{enumerate}
\item Every even Dirac structure on $M$ is the union of a presymplectic manifold and a foliated Poisson manifold. These manifolds are glued along the region $M_{(1, 0)}$.
\item Every odd Dirac structure on $M$ is the union of a Poisson manifold and a foliated presymplectic manifold. These manifolds are glued along the region $M_{(0, 1)}$.
\end{enumerate}
\end{theorem}


\begin{thebibliography}{100}

\bibitem{C} Courant, T. \textit{Dirac manifolds.} Trans. Amer. Math. Soc. 319, no. 2, pp. 631-661, 1990.

\bibitem{CW} Courant, T. and Weinstein, A. \textit{Beyond Poisson Structures} Action hamiltoniennes de groupes. Troisieme theoreme de Lie (Lyon, 1986), pp.39--49, Travaux en Cours, 27, Hermann, Paris, 1988.

\bibitem{G} Gualtieri, M. \textit{Generalized Complex Geometry} Ann. of Math. 174 (2011), pp. 75-123
\end{thebibliography}
\end{document}